\documentclass[12pt]{article}

\usepackage{amsmath}
\usepackage{amsthm}
\usepackage{amsfonts}
\usepackage{amssymb}

\numberwithin{equation}{section}

\newtheorem{thm}{Theorem}[section]
\newtheorem{assertion}[thm]{Proposition}
\newtheorem{lemma}[thm]{Lemma}
\newtheorem{rem}[thm]{Remark}

\theoremstyle{remark}

\newcommand{\suin}{\underset{n\geq 0}{\sup}}

\newcommand{\zi}{Z_{\infty}}

\newcommand{\lin}{\underset{n\rightarrow\infty}{\lim}}

\newcommand{\lix}{\underset{x\rightarrow\infty}{\lim}}

\newcommand{\mmp}{\mathbb{P}}
\newcommand{\me}{\mathbb{E}}
\newcommand{\mn}{\mathbb{N}}

\newcommand{\mr}{\mathbb{R}}

\newcommand{\mm}{\mathcal{Z}}

\def\gl{\buildrel \rm def\over =}
\def\eqdist{\ {\buildrel d\over =}\ }
\def\1{{\bf 1}}

\newcommand{\bfP}{\mathbf{P}}
\newcommand{\bfQ}{\mathbf{Q}}

\newcommand{\whL}{\widehat{L}}
\newcommand{\whW}{\widehat{W}}
\newcommand{\whS}{\widehat{S}}
\newcommand{\whbfL}{\widehat{\mathbf{L}}}
\newcommand{\whbfT}{\widehat{\mathbf{T}}}
\newcommand{\whcalN}{\widehat{\mathcal{N}}}

\begin{document}

%
%

\title{\Large\bf A Log-Type Moment Result for Perpetuities and Its Application
to Martingales in Supercritical Branching Random Walks}\date{}
\author{Gerold Alsmeyer\footnote{email address:
gerolda@math.uni-muenster.de}\\\small{\emph{Institut f\"{u}r Mathematische
Statistik},
\emph{Westf\"{a}lische Wilhelms-Universit\"{a}t M\"{u}nster}},\vspace{-.2cm}\noindent\\
\small{\emph{Einsteinstrasse 62, 48149 M\"{u}nster, Germany}}\vspace{-.3cm}\noindent\\ \\
Alexander
Iksanov\footnote {email address: iksan@unicyb.kiev.ua}\\\small{\emph{Faculty of
Cybernetics},
\emph{National T. Shevchenko University}},\vspace{-.2cm}\noindent\\
\small{\emph{01033 Kiev, Ukraine}}} \maketitle
\begin{quote}\vspace{-.3cm}\noindent
\noindent {\footnotesize \textbf{SUMMARY.}} Infinite sums of i.i.d.\ random
variables discounted by a multiplicative random walk are called perpetuities
and have been studied by many authors. The present paper provides a log-type
moment result for such random variables under minimal conditions which is then
utilized for the study of related moments of a.s.\ limits of certain martingales
associated with the supercritical branching random walk. The connection, first observed
by the second author in \cite{Iks041}, arises upon consideration of a size-biased
version of the branching random walk originally introduced by Lyons \cite{Lyons}. We
also provide a necessary and sufficient condition for uniform integrability of these
martingales in the most general situation which particularly means that the classical
(LlogL)-condition is not always needed.
\end{quote}

\section{Introduction and results}

The principal purpose of this article is to provide a log-type moment result
for the limit of iterated i.i.d.\ random linear functions, called {\it
perpetuties}. It is given as Theorem \ref{perp} in the following subsection along
with all necessary facts about the model. A similar result (Theorem
\ref{brw_conc2}) will then be formulated for the a.s.\ limit of a well-known martingale
associated with the branching random walk introduced in Subsection 1.2.
As will be explained in Section \ref{ui}, the connection between these at first glance
unrelated models pops up when studying the weighted random tree associated with
the branching random walk under the so-called size-biased measure. It does not
take by surprise that this connection, once established, can be utilized to
obtain moment results in the branching model by resorting to corresponding
ones for perpetuities.

\subsection{Perpetuities}

Given a sequence $\{(M_n, Q_n): n=1,2,...\}$ of i.i.d.\ $\Bbb{R}^{2}$-valued
random vectors with generic copy $(M, Q)$, put
$$\Pi_0\ \gl\ 1\quad\text{and}\quad\Pi_n\ \gl\ M_1 M_2 \cdots M_n,
\quad n=1,2,... $$
and
$$ Z_n\ \gl\ \sum_{k=1}^n \Pi_{k-1}Q_k,\quad n=1,2,... $$
The random discounted sum
\begin{equation}\label{per}
\zi\ \gl\ \sum_{k\ge 1} \Pi_{k-1}Q_k,
\end{equation}
obtained as the a.s.\ limit of $Z_n$
under appropriate conditions (see Proposition \ref{exper} below), is called
perpetuity and of interest in various fields of applied probability like insurance and
finance, the study of shot-noise processes or, as will be seen further on, of
branching random walks. The law of $\zi$ appears also quite naturally as the stationary
distribution of the (forward) iterated function system
\begin{equation*} 
\Phi_{n}\ \gl\ \Psi_{n}(\Phi_{n-1})\ =\ \Psi_{n}\circ...\circ\Psi_{1}
(\Phi_{0}),\quad n=1,2...,
\end{equation*}
where $\Psi_{n}(t)\gl Q_{n}+M_{n}t$ for $n=1,2,...$ and $\Phi_{0}$ is independent
of $\{(M_n, Q_n): n=1,2,...\}$. Due to the recursive structure of this Markov chain,
it forms solution of the stochastic fixed point equation
\begin{equation*} 
\Phi\ \eqdist\ Q+M\Phi
\end{equation*}
where as usual the variable $\Phi$ is assumed to be independent of $(M,Q)$.
Let us finally note that $\zi$ may indeed be obtained as the
a.s.\ limit of the associated backward system when started at $\Phi_{0}\equiv 0$, i.e.
\begin{equation*} 
\zi\ =\ \lim_{n\to\infty}\Psi_{0}\circ...\circ\Psi_{n}(0).
\end{equation*}

Goldie and Maller \cite{GolMal} gave the following complete characterization of
the a.s.\ convergence of the series in (\ref{per}). For $x>0$, define
\begin{equation}\label{A(x)}
A(x)\ \gl\ \int_0^x \mmp\{-\log |M|>y\}\ dy\ =\ \me\min\big(\log^{-}|M|,x\big)
\end{equation}
and then $J(x)\gl x/A(x)$. In order to have $J(x)$ defined on the whole real line,  put
$J(x)\gl 0$ for $x<0$ and $J(0)\gl\lim_{x\downarrow 0}J(x)=1/\mmp\{|M|<1\}$.

\begin{assertion}\label{exper} {\rm (\cite{GolMal}, Theorem 2.1)}
Suppose
\begin{equation}\label{nonzero}
\mmp\{M=0\}=0\quad\text{and}\quad\mmp\{Q=0\}<1.
\end{equation}
Then
\begin{equation}\label{cond2000}
\lin \Pi_n\,=\,0\ \text{a.s.}\quad\text{and}\quad\me J\big(\log^{+}|Q|\big)
\,<\,\infty,
\end{equation}
and
\begin{equation}\label{conv}
\zi^{*}\ \gl\ \sum_{n\ge 1}|\Pi_{n-1}Q_n|\ <\ \infty\quad\text{a.s.}
\end{equation}
are equivalent conditions, and they imply
\begin{equation*}
\lin Z_n\,=\,\zi\text{ a.s.\quad and}\quad|\zi|\ <\ \infty\quad\text{a.s.}
\end{equation*}
Moreover, if
\begin{equation}\label{degen}
\mmp\{Q+Mc=c\}<1\quad\text{for all } c\in \mr,
\end{equation}
and if at least one of the conditions in (\ref{cond2000}) fails to hold, then
$\lin |Z_n|=\infty$ in probability.
\end{assertion}

Condition (\ref{cond2000}) holds particularly true if
\begin{equation}\label{simple1}
\me \log |M|\in (-\infty, 0)\quad\text{and}\quad\me\log^+|Q|<\infty,
\end{equation}
and for this special case results on the finiteness of certain log-type
moments of $\zi$ were derived in \cite{Iks06} and \cite{IksRos}.
To extend those results to the general situation with (\ref{nonzero}) being the only
basic assumption is one purpose of the present paper.


Let the function $b:\mr^+\to\mr^+$ be measurable, locally bounded and regularly
varying at $\infty$ with exponent $\alpha>0$. Functions $b$ of interest
in the following result are, for instance, $b(x)=x^{\alpha}\log_{k}x$ or
$b(x)=x^{\alpha}\exp(\beta\log^{\gamma}x)$ for $\beta\ge 0$, $0<\gamma<1$
and $k\in\Bbb{N}$, where $\log_{k}$ denotes $k$-fold iteration of the
logarithm.
\begin{thm}\label{perp} Suppose (\ref{nonzero}). Then $\lim_{n\to\infty}\Pi_n=0$ a.s.,
\begin{equation}\label{M}
\me b\big(\log^{+}|M|\big)J\big(\log^{+}|M|\big)<\infty
\end{equation}
and
\begin{equation}\label{Q}
\me b\big(\log^{+}|Q|\big)J\big(\log^{+}|Q|\big)<\infty
\end{equation}
together imply
\begin{equation}\label{Z}
\me b(\log^+|\zi|)<\infty.
\end{equation}
Conversely, if $\zi$ is a.s.\ finite and nondegenerate, then (\ref{Z}) implies
(\ref{M}) and (\ref{Q}).
\end{thm}
Replacing $\lim_{n\to\infty}\Pi_{n} =0$ a.s.\ with the stronger condition
$\me\log|M|\in (-\infty, 0)$, this result is stated as Theorem 3 in
\cite{Iks06}, and our proof also fixes a minor flaw appearing in the proof given
there.

Since (\ref{M}) and (\ref{Q}) are conditions in terms of the absolute values of $M$ and
$Q$, the first conclusion of Theorem \ref{perp} remains valid when replacing (\ref{Z})
with the stronger assertion
\begin{equation}\label{Z*}
\me b(\log^+\zi^{*})<\infty.
\end{equation}
If $\Pi_{n}\to 0$ a.s.\ and if $\zi$ and $\zi^{*}$ are both a.s.\ finite
and nondegenerate, this leads us to the conclusion that (\ref{Z}) and (\ref{Z*}) are
actually equivalent. A similar conclusion has been obtained in \cite{AlIkRo} for the
case of ordinary moments (viz.\ $b(\log x)=x^{p}$ for some $p>0$), see
Theorem 1.4 there.

\subsection{The branching random walk and its intrinsic martingales}\label{brw}

In the following we give a short description of the standard branching random
walk, its intrinsic martingales and an associated multiplicative random walk.

Consider a population starting from one ancestor located at the origin
and evolving like a Galton-Watson process but with the generalization that
individuals may have infinitely many children. All individuals are residing in
points on the real line, and the displacements of children relative to their
mother are described by a point process $\mm=\sum_{i=1}^{N} \delta_{X_i}$ on
$\mr$. Thus $N\gl\mm(\mr)$ gives the total number of offspring of the considered
mother and $X_{i}$ the displacement of the $i$-th child. The displacement
processes of all population members are supposed to be independent copies of
$\mm$. We further assume $\mm(\{-\infty\})=0$ and $\me N>1$ (supercriticality)
including the possibility $\mmp\{N=\infty\}>0$ as already stated above. If
$\mmp\{N<\infty\} =1$, then the population size process forms an ordinary
Galton-Watson process. Supercriticality ensures survival of the population with
positive probability.

For $n=0,1,...$ let $\mm_n$ be the point process
that defines the positions on $\mr$ of the individuals of the $n$-th generation,
their total number given by $\mm_{n}(\mr)$.
The sequence $\{\mm_n: n=0,1,...\}$ is called
\emph{branching random walk} (BRW).

Let $\mathbf{V}\gl\bigcup_{n=0}^\infty \mn^n$ be the infinite Ulam-Harris tree
of all finite sequences $v=v_1...v_n$ (shorthand for $(v_{1},...,v_{n})$),
with root $\varnothing$ ($\mn^0\gl\{\varnothing\}$) and edges connecting
each $v\in\mathbf{V}$ with its successors $vi$, $i=1,2,...$ The length of $v$ is
denoted as $|v|$. Call $v$ an individual and $|v|$ its generation number.
A BRW $\{\mm_n: n=0,1,...\}$ may now be represented as a random labeled subtree
of $\mathbf{V}$ with the same root. This subtree $\mathbf{T}$ is obtained recursively
as follows: For any $v\in\mathbf{T}$, let $N(v)$ be the number of its
successors (children) and $\mm(v)\gl\sum_{i=1}^{N(v)}\delta_{X_{i}(v)}$ denote
the point process describing the displacements of the children $vi$ of $v$
relative to their mother. By assumption, the $\mm(v)$ are independent copies
of $\mm$. The Galton-Watson tree associated with this model is now given by
$$ \mathbf{T}\gl\{\varnothing\}\cup\{v\in\mathbf{V}\backslash\{\varnothing\}:
v_{i}\le N(v_{1}...v_{i-1})\text{ for }i=1,...,|v|\}, $$
and $X_{i}(v)$ denotes the label attached to the edge $(v,vi)\in\mathbf{T}\times
\mathbf{T}$ and describe the displacement of $vi$ relative to $v$. Let us stipulate
hereafter that $\sum_{|v|=n}$ means summation over all vertices of $\mathbf{T}$ (not
$\mathbf{V}$) of length $n$. For
$v=v_{1}...v_{n}\in\mathbf{T}$, put $S(v)\gl\sum_{i=1}
^{n}X_{v_{i}}(v_{1}...v_{i-1})$. Then $S(v)$ gives the position of $v$ on the
real line (of course, $S(\varnothing)=0$), and
$\mm_{n}=\sum_{|v|=n}\delta_{S(v)}$ for all $n=0,1,..$.

Suppose there exists $\gamma>0$ such that
\begin{equation}\label{gam2000}
m(\gamma)\ \gl\ \me \int_\mr e^{\gamma x}\,\mm(dx)\in (0,\infty).
\end{equation}
For $n=1,2,...$, define $\mathcal{F}_n\gl\sigma(\mm(v):|v|
\le n-1)$, and let $\mathcal{F}_0$ be the trivial
$\sigma$-field. Put
\begin{equation}\label{Wn2000}
W_n\ \gl\ m(\gamma)^{-n}\int_\mr e^{\gamma x}\,\mm_n(dx)
\ =\ m(\gamma)^{-n}\sum_{|v|=n}e^{\gamma S(v)}\ =\ \sum_{|v|=n}L(v),
\end{equation}
where $L(v)\gl e^{\gamma S(v)}/m(\gamma)^{|v|}$. Notice that the dependence of
$W_{n}$ on $\gamma$\nobreak\ has been suppressed. The sequence
$\{(W_n,\mathcal{F}_n):n=0,1,...\}$ forms a non-negative martingale with mean one
and is thus a.s.\ convergent with limiting variable $W$, say, satisfying $\me W
\le 1$. It has been extensively studied in the literature, but first results were
obtained in \cite{King} and \cite{Big1}. Note that
$\mmp\{W>0\}>0$ if, and only if, $\{W_n:n=0,1,...\}$ is
uniformly integrable. While uniform integrability is clearly sufficient, the
necessity hinges on the well known fact that $W$
satisfies the stochastic fixed point equation
\begin{equation}\label{fixp}
W\ =\ \sum_{|v|=n}L(v)W(v)\quad\hbox{a.s.}
\end{equation}
for $n=1,2,...$, where the $W(v)$, $|v|=n$, are i.i.d.\ copies of $W$ that
are also independent of $\{L(v):|v|=n\}$, see e.g.\ \cite{BiggKypr97}. In fact
$W(v)$ is nothing but the a.s.\ limit of the martingale $\{\sum_{|w|=m}{L(vw)\over
L(v)}:m=0,1,...\}$ which forms the counterpart of $\{W_{n}:n=0,1,...\}$, but for
the subtree of $\mathbf{T}$ rooted at $v$.

Our goal is to study certain moments of $W$ in the nontrivial situation where
$\{W_n:n=0,1,...\}$ is uniformly integrable. For the latter to hold, Theorem
\ref{brwL1} below provides us with a necessary and sufficient condition, again
under no additional assumptions on the BRW beyond (\ref{gam2000}). In order to
formulate it, we first need to introduce a multiplicative random walk associated
with our model. Let
$M$ be a random variable with distribution defined by
\begin{equation}\label{Zm}
\mmp\{M\in B\}\ \gl\ \me\left[\sum_{|v|=1}L(v)\delta_{L(v)}(B)\right],
\end{equation}
for any Borel subset $B$ of $\mr^+$. Notice that the right-hand side of (\ref{Zm})
does indeed define a probability distribution because $\me\sum_{|v|=1}L(v)=\me
W_{1}=1$. More generally, we have (see e.g.\ \cite{BiggKypr97}, Lemma 4.1)
\begin{equation}
\label{Z2}
\mmp\{\Pi_{n}\in B\}\ =\ \me\left[\sum_{|v|=n}L(v)\delta_{L(v)}(B)
\right],
\end{equation}
for each $n=1,2,...$, whenever $\{M_k: k=1,2,...\}$ is a family of independent
copies of $M$ and $\Pi_{n}\gl\prod_{k=1}^{n}M_{k}$. It is important to note
that
\begin{equation}\label{Z1}
\mmp\{M=0\}=0\quad\text{and}\quad\mmp\{M=1\}<1.
\end{equation}
The first assertion follows since, by (\ref{Zm}), $\mmp\{M>0\}=\me W_{1}=1$. As
\nobreak\ for the second, observe that $\mmp\{M=1\}=1$ implies
$\me\sum_{|v|=1}L(v)\1_{\{L(v)\ne 1\}}=0$ which in combination with $\me
W_{1}=1$ entails that the point process $\mm$ consists of only one point $u$ with
$L(u)=1$. This contradicts the assumed supercriticality of the BRW.

Not surprisingly, the chosen notation for the multiplicative random walk
associated with the given BRW as opposed to the notation in the previous
subsection is intentional, and we also keep the definitions of $J(x)$ and $A(x)$
from there, see (\ref{A(x)}) and thereafter.

\begin{thm}\label{brwL1} The martingale $\{W_n:n=0,1,...\}$ is uniformly
integrable if, and only if, the following two conditions hold true:
\begin{equation} \label{Z3}
\lim_{n\to\infty}\Pi_{n}=0\quad\text{a.s.}
\end{equation}
and
\begin{equation}\label{mea}
\me W_{1}J(\log^{+}W_{1})\ =\ \int_{(1,\infty)}xJ(\log x)
\ \mmp\{W_{1}\in dx\}\ <\ \infty.
\end{equation}
\end{thm}

There are three distinct cases in which conditions (\ref{Z3}) and (\ref{mea})
hold simultaneously:
\newline (A1)\quad
$\me \log M \in (-\infty, 0)$ and $\me W_1 \log^+ W_1<\infty$;
\newline (A2)\quad $\me \log M=-\infty$ and $\me W_{1}J
(\log^{+}W_{1})<\infty$;
\newline (A3)\quad $\me \log ^+M=\me \log^- M=+\infty$,
$\me W_{1}J(\log^{+}W_{1})<\infty$, and
$$ \me J\big(\log^{+}M\big)\ =\
\int_{(1,\infty )}\dfrac{\log x}{\int_0^{\log x}
\mmp\{-\log M>y\}\,dy}\ \mmp\{M \in dx\}\ <\ \infty. $$
For the case (A1), Theorem \ref{brwL1} is due to
Biggins \cite{Big1} and Lyons \cite{Lyons}, see also \cite{Kuhl2}.
In the present form, the result has also been stated as Proposition 1 in
\cite{IksRos} (with a minor misprint), however without proof and a reference
to the proof of Theorem 2 in \cite{Iks041}
instead. But the latter theorem was formulated\nobreak\ in terms of fixed
points rather than martingale convergence which somewhat obscures how to extract
the necessary arguments. On the other hand, the study of uniform integrability
has a long history, going back to the famous Kesten-Stigum theorem
\cite{KestStig} for ordinary Galton-Watson processes and the pioneering work by
Biggins \cite{Big1} for the BRW, and followed later by work in \cite{Liu} and
\cite{Lyons}. We have therefore decided to include a complete (and rather short)
proof here.

The existence of moments of $W$ was studied in quite a number of articles, see
\cite{AlsKu},\cite{Big1},\cite{BiDo75},\cite{Iks06},\cite{IksRos},\cite{
Liu2000},\cite{RosVat}. The following theorem, which is our second main
moment-type result, goes beyond the afore-mentioned ones in that it does not
restrict to case (A1) of Theorem \ref{brwL1}. The function
$b(x)$ occurring here is of the type stated before Theorem \ref{perp}.

\begin{thm}\label{brw_conc2}
If $\lim_{n\to\infty}\Pi_{n}=0$ a.s. and
\begin{equation}\label{ZZZ}
\me W_{1}b\big(\log^{+}W_{1}\big)J(\log^{+}W_{1})\ <\ \infty,
\end{equation}
then $\{W_n:n=0,1,...\}$ is uniformly integrable and
\begin{equation}\label{ZZZZ}
\me Wb(\log^+ W)<\infty.
\end{equation}
Conversely, if (\ref{ZZZZ}) holds and $\mmp\{W_1=1\}<1$, then (\ref{ZZZ}) holds.
\end{thm}

An interesting aspect of this theorem is that it provides conditions
for the existence of $\Phi$-moments of $W$ for $\Phi$ slightly beyond
$\mathcal{L}_{1}$ without assuming the (LlogL)-condition to ensure uniform
integrability. The latter condition is a standing assumption in a related article by
Alsmeyer and Kuhlbusch \cite{AlsKu} where a similar but more general result (as
for the functions $\Phi$) is proved, see Theorem 1.2 there.

There are basically two probabilistic approaches towards finding
conditions for the existence of $\me \Phi(W)$ for suitable functions
$\Phi$. The method of this paper, worked out in \cite{Iks041} and
\cite{IksRos}, hinges on getting first a moment-type result for perpetuities
(here Theorem \ref{perp}) and then translating it into the framework of branching
random walks. This is accomplished by an appropriate change of measure argument
(see the proof of Theorem \ref{brwL1}). The second approach, first used in
\cite{AlsRos} for Galton-Watson processes and further elaborated in \cite{AlsKu},
relies on the observation that BRW's bear a certain double martingale structure
which allows the repeated application of the convex function inequalities due
to Burkholder, Davis and Gundy (see e.g.\ \cite{Chow}) for
martingales. Both approaches have their merits and limitations. Roughly speaking,
the double martingale argument requires as indispensable ingredients only that
$\Phi$ be convex and at most of polynomial growth. On the other hand, it also
comes with a number of tedious technicalities caused by the repeated application
of the convex function inequalities. The basic tool of the method used here is
only Jensen's inequality for conditional expectations, but it relies heavily on
the existence of a nonnegative concave function $\Psi$ that is equivalent at
$\infty$ to the function $\Phi(x)/x$. This clearly imposes a strong restriction on
the growth of $\Phi$.

The rest of the paper is organized as follows. Section \ref{func} collects the
relevant properties of the functions involved in our analysis, notably $b(x)$,
$b(\log x)$ and $A(x)$, followed in Section \ref{aux} by some preliminary work
needed for the proofs of Theorems \ref{perp} and \ref{brw_conc2}. In
particular, a number of moment results for certain functionals of multiplicative
random walks are given there which may be of independent interest (see Lemma
\ref{collection}). Theorem \ref{perp} is proved in Section
\ref{pro1}, while Section \ref{ui} contains the proofs of Theorems \ref{brwL1}
and \ref{brw_conc2}.

\section{Properties of the functions involved}\label{func}

In this section, we gather some relevant properties of the functions $b(x)$,
$A(x)$ and $J(x)=x/A(x)$ needed in later on. Recall from (\ref{A(x)})
the definition of $A(x)$ and that $b:\mr^+\to\mr^+$ is measurable, locally bounded
and regularly varying at $\infty$ with exponent $\alpha>0$ and thus of the
form $b(x)=x^{\alpha}\ell(x)$ for some slowly varying function $\ell(x)$. By the
Smooth Variation Theorem (see Thm.\ 1.8.2 in \cite{BGT}), we may assume without
loss of generality that $b(x)$ is smooth with $n$th derivative $b^{(n)}(x)$
satisfying
$$ x^{n}b^{(n)}(x)\sim \alpha(\alpha-1)\cdot...\cdot (\alpha-n+1) b(x) $$
for all $n\ge 1$, where $f\sim g$ has the usual meaning that $\lim_{x\to\infty}
f(x)/g(x)=1$. By Lemma 1 in \cite{Als}, $b(x)$ may further be chosen in such a way
that
\begin{equation}\label{subadd0}
b(x+y)\ \le\ C\big(b(x)+b(y)\big)
\end{equation}
for all $x,y\in\Bbb{R}^{+}$ and some $C\in (0,\infty)$. The smoothness of $b(x)$
(and thus of $\ell(x)$) and property (\ref{subadd0}) will be standing assumptions
throughout without further notice.

Before giving a number of lemmata, let us note the obvious facts that
\medskip
\newline (P1)\quad $A(x)$ is nondecreasing,
\newline (P2)\quad $J(x)$ is nondecreasing with $\lim_{x\to\infty}J(x)=\infty$,
and
\newline (P3)\quad $J(x)\sim J(x+a)$ for any fixed $a>0$.

\begin{lemma}\label{function} There exist smooth nondecreasing and concave
functions $f$ and $g$ on $\Bbb{R}^{+}$ with $f(0)=g(0)=0$,
$\lim_{x\to\infty} f(x)=\lim_{x\to\infty} g(x)=\infty$,
$f^\prime(0+)<\infty$ and $g'(0+)<\infty$ such that $b(\log
x)\sim f(x)$ and $b(\log x)\log x \sim g(x)$. Moreover,
\begin{equation}\label{ineq}
f(xy)\le C(f(x)+f(y))
\end{equation}
for all $x,y\in\Bbb{R}^{+}$ and some $C\in (0,\infty)$.
\end{lemma}

\begin{proof} For each $c>0$, we have that $\Lambda_{c}(x)\gl b(\log(c+x))
-b(\log c)$ satisfies $\Lambda_{c}(0)=0$, $\Lambda_{c}(x)\sim b(\log x)$ and
$\Lambda_{c}'(x)={b'(\log(c+x))\over c+x}\sim {\alpha b(\log(c+x))\over (c+x)
\log(c+x)}$. We thus see that $\Lambda_{c}'(x)$ is regularly varying of order $-1$ and,
for $c$ sufficiently large, nonincreasing on
$\Bbb{R}^{+}$ with $\Lambda_{c}'(0+)=c^{-1}b'(\log c)\in (0,\infty)$. Similar
statements hold true for $\Lambda_{c}(x)\log(c+x)\sim b(\log x)\log x$.
Since $\Lambda_{c}(e^{x})\sim b(x)$ and $b(x)$ satisfies (\ref{subadd0}), it
is readily verified that $\Lambda_{c}(x)$ satisfies (\ref{ineq}).
Consequently, the lemma follows upon choosing $f(x)=\Lambda_{c}(x)$ and
$g(x)=\Lambda_{c}(x)\log(c+x)$ for sufficiently large $c$.
\end{proof}

\begin{lemma}\label{subadd}
Let $g$ be as in Lemma \ref{function}. Then $\phi(x)\gl g(x)/A(\log
(x+1))$ is subadditive on $\Bbb{R}^{+}$, i.e.\ $\phi(x+y)\le \phi(x)+\phi(y)$ for
all $x,y\geq 0$, and $f(x)J(\log x)\sim\phi(x)$.
\end{lemma}

\begin{proof} Since $g$ is concave, $g(\alpha x)\geq \alpha g(x)$ for each
$\alpha\in (0,1)$ and $x\geq 0$. Hence we infer with the help of (P1)
\begin{equation}\label{star}
\phi(\alpha x)\geq \alpha \phi(x) \ \text{for every} \ \alpha\in (0,1) \
\text{and} \ x\geq 0
\end{equation} which implies subadditivity via $\phi(x)+\phi(y)\ge
[{x\over x+y}+{y\over x+y}]\phi(x+y)=\phi(x+y)$. The asymptotic result
follows from $g(x)\sim f(x)\log x\sim f(x)\log(x+1)$ (see Lemma \ref{function})
which implies
$$ \phi(x)\ \sim\ f(x)J(\log(x+1))\ \sim\ f(x)J(\log x) $$
having utilized (P2) and (P3) for the last asymptotic equivalence.
\end{proof}

\begin{lemma}\label{slowvar}
The function $\phi$ in Lemma \ref{subadd} is slowly varying at $\infty$ and
satisfies $\phi(x)\sim \phi(x+b)$ for any fixed $b\in \mr$. Furthermore,
\begin{equation}\label{submulti}
\phi(xy)\ \le\ C(\phi(x)+\phi(y))
\end{equation}
for all $x,y\in\mr^{+}$ and a suitable constant $C\in (0,\infty)$.
\end{lemma}

\begin{proof} We must check $\lix \phi(xy)/\phi(x)=1$ for $y>1$.
By the previous lemma, we have
$$ \dfrac{\phi(xy)}{\phi(x)}\ \sim\ \dfrac{f(xy)}{f(x)}\,\dfrac{J(\log x+\log y)}
{J(\log x)}, $$
which yields the desired conclusion because $f(x)\sim b(\log x)$ is slowly varying
and, by (P3), $J(\log x+\log y)\sim J(\log x)$ for any fixed $y$. The second
assertion follows as a simple consequence so that we turn directly to
(\ref{submulti}). Fix $K\in\Bbb{N}$ so large that
${\phi(x)\over f(x)J(\log x)}\in [1/2,2]$ for all $x\ge K$ and use the subadditivity of
$\phi$ to infer in the case $x\wedge y\le K$
\begin{equation}\label{eq1}
\phi(xy)\ \le\ \phi(K(x\vee y))\ \le K(\phi(x)\vee\phi(y))
\ \le\ K(\phi(x)+\phi(y)).
\end{equation}
Note next that $J$ as a nondecreasing sublinear function satisfies
$J(x+y)\le C(J(x)+J(y))$ for all $x,y\in\mr^{+}$. By combining this with the
monotonicity of $f,J$ and inequality (\ref{ineq}), we obtain if $x>K$ and $y>K$
(thus $xy>K$)
\begin{eqnarray}\label{eq2}
\phi(xy)&\le&2f(xy)J(\log x+\log y)\nonumber\\
&\le&2C(f(x)+f(y))(J(\log x)+J(\log y))\nonumber\\
&\le&8C(f(x)J(\log x)\vee f(y)J(\log y))\nonumber\\
&\le&16C(\phi(x)+\phi(y)),
\end{eqnarray}
for a suitable constant $C\in (0,\infty)$. A combination if (\ref{eq1}) and
(\ref{eq2}) yields (\ref{submulti}) (with a suitable $C$).
\end{proof}

\section{Auxiliary results}\label{aux}


In the notation of Subsection 1.1 and always assuming (\ref{nonzero}), let us
consider the situation where $|\zi|<\infty$ and the nondegeneracy condition
(\ref{degen}) is in force. Then $\lim_{n\to\infty}\Pi_{n}=0$ by Proposition
\ref{exper}, and
\begin{equation}\label{pereqdistr}
\zi\ =\ Q_{1}+M_{1}\zi^{(1)}\ =\ Q^{(m)}+\Pi_{m}\zi^{(m)},
\end{equation}
holds true for each $m\ge 1$, where (setting $\Pi_{k:l}\gl M_{k}\cdot...\cdot
M_{l}$)
\begin{equation}\label{iterate}
Q^{(m)}\ \gl\ \sum_{k=1}^{m}\Pi_{k-1}Q_{k}\quad\text{and}\quad
\zi^{(m)}\ \gl\ Q_{m+1}+\sum_{k\ge m+2}\Pi_{m+1:k-1}Q_k.
\end{equation}
Here $\zi^{(m)}$ constitutes a copy of $\zi$ independent of $(M_1,Q_1),...,
(M_{m},Q_{m})$. We thus see that $\zi$ may also be viewed as the perpetuity
generated by i.i.d.\ copies of $(\Pi_{m},Q^{(m)})$ for any fixed $m\ge 1$.
We may further replace $m$ by any a.s. finite stopping time $\sigma$ to obtain
\begin{equation}\label{genfpeq}
\zi\ =\ \sum_{k=1}^{\sigma}\Pi_{k-1}Q_{k}\ +\ \Pi_{\sigma}\zi^{(\sigma)},
\end{equation}
where $Q^{(\sigma)}\gl\sum_{k=1}^{\sigma}\Pi_{k-1}Q_{k}$ and $\zi^{(\sigma)}$ is a
copy of $\zi$ independent of $\sigma$ and $\{(M_{n},Q_{n}): 1\le n\le\sigma\}$
(and thus of $(\Pi_{\sigma},Q^{(\sigma)})$). For our purposes, a relevant choice
of $\sigma$ will be
\begin{equation}\label{ladder}
\sigma\ \gl\ \inf\{n\ge 1:|\Pi_{n}|\le 1\},
\end{equation}
which is nothing but the first (weakly) ascending ladder epoch for the random walk
$S_{n}\gl -\log|\Pi_{n}|$, $n=0,1,...$

\begin{lemma}\label{perpsup} Let $\zi$ be nondegenerate and $f$ be a function as in
Lemma \ref{function}. Define
$$ Q_{n}^{(2)}\ \gl\ Q_{2n-1}+M_{2n-1}Q_{2n} $$
for $n\ge 1$ and
let $\overline{Q}_{n}^{(2)}$ be a conditional symmetrization of $Q_{n}^{(2)}$
given $M_{2n-1}M_{2n}$. Then $\me f(|\zi|)<\infty$ implies
\begin{eqnarray}
&&\me f(|Q|)<\infty\quad\text{and}\quad\me f(|M|)<\infty,\label{cond2001}\\
&&\me f\Big(\sup_{n\ge 1}|\Pi_{n-1}Q_{n}|\Big)<\infty,
\label{cond2002}\\
&&\me f\Big(\sup_{n\ge 1}|\Pi_{2n-2}\overline{Q}_{n}^{(2)}|\Big)<\infty,
\label{cond3001}\\
&&\me f\Big(\sup_{n\ge 0}|\Pi_n|\Big)<\infty.\label{cond3002}
\end{eqnarray}
\end{lemma}
\begin{proof} It has been shown in \cite{AlIkRo} that, under the given
assumptions, the distribution of $\overline{Q}_{n}^{(2)}$ is nondegenerate,
\begin{equation}\label{symm}
\mmp\Big\{\sup_{k\ge 1}|\Pi_{2k-2}\overline{Q}_{k}^{(2)}|>x\Big\}
\ \le\ 4\,\mmp\{|\zi|>x/2\}
\end{equation}
for all $x>0$ (see (28) there) and
\begin{equation}\label{tailin}
\mmp\Big\{\sup_{k\ge 0}|\Pi_{2k}|>x\Big\}\ \le\ 2\,\mmp\Big\{
\sup_{k\ge 1}|\Pi_{2k-2}\overline{Q}_{k}^{(2)}|>cx\Big\}
\end{equation}
for all $x>0$ and a suitable $c\in (0,1)$ (see Lemma 2.1 of \cite{AlIkRo}). By our
standing assumption (\ref{nonzero}), we can choose
$0<\rho<1$ so small that $\kappa\gl\mmp\{|M|>\rho\}>0$. With the help
of the above tail inequalities we now infer (\ref{cond3001}) and thereupon
(\ref{cond3002}) because
\begin{eqnarray*}
\mmp\Big\{\sup_{k\ge 0}|\Pi_{2k}|>\rho x\Big\} &\ge& \mmp\Big\{\sup_{k\ge 1}
|\Pi_{2k}|>\rho x,|M_{1}|>\rho\Big\}\\
&\ge& \mmp\Big\{\sup_{k\ge 1}|\Pi_{2:2k}|>x,|M_{1}|>\rho\Big\}\\
&=& \kappa\,\mmp\Big\{\sup_{k\ge 1}|\Pi_{2k-1}|>x\Big\}
\end{eqnarray*}
and thus
\begin{eqnarray*}
\mmp\Big\{\sup_{k\ge 0}|\Pi_{k}|>2x\Big\}&\le&\mmp\Big\{\sup_{k\ge
0}|\Pi_{2k}|>x\Big\}+\mmp\Big\{\sup_{k\ge 1}|\Pi_{2k-1}|>x\Big\}\\
&\le&(1+\kappa^{-1})\mmp\Big\{\sup_{k\ge 0}|\Pi_{2k}|>\rho x\Big\}
\end{eqnarray*}
for all $x>0$. Next, $\me f(|M|)<\infty$ follows from (\ref{cond3002}) and
$|M_{1}|\le\sup_{n\ge 0}|\Pi_{n}|$. As for $\me f(|Q|)<\infty$, we recall
from (\ref{pereqdistr}) that $\zi=Q_{1}+M_{1}\zi^{(1)}$. Hence
$$ \me f(|Q_{1}|)\ \le\ \me f(|\zi|)+\me f(|M_{1}\zi^{(1)}|)\ \le\ C\Big(\me
f(|\zi|)+\me f(|M_{1}|)\Big)\ <\ \infty $$
for a suitable $C\in (0,\infty)$, where subadditivity of $f$ has been used for the
first inequality and (\ref{ineq}) for the second one.

Finally, we must verify (\ref{cond2002}). With $m_{0}$ denoting a median of $\zi$,
Goldie and Maller (see \cite{GolMal}, p. 1210) showed that
$$ \mmp\Big\{\sup_{n\ge 1}|Z_{n}+\Pi_{n}m_{0}|>x\Big\}\ \le\ 2\,\mmp\{|\zi|
\ge x\} $$
for all $x>0$. Hence $\me f(\sup_{n\ge 1}|Z_{n}+\Pi_{n}m_{0}|)
\le 2\,\me f(|\zi|)<\infty$. Now
$$ \Pi_{n-1}Q_{n}\ =\ (Z_{n}+\Pi_{n}m_{0})-(Z_{n-1}-\Pi_{n-1}m_{0})
+m_{0}(\Pi_{n-1}-\Pi_{n}) $$
implies (as $Z_{0}=0$ and $\Pi_{0}=1$)
$$ \sup_{n\ge 1}|\Pi_{n-1}Q_{n}|\ \le\ 2\Big(\sup_{n\ge 0}|Z_{n}+\Pi_{n}m_{0}|+
|m_{0}|\sup_{n\ge 0}|\Pi_{n}|\Big)+|m_{0}|, $$
and this gives the desired conclusion by (\ref{cond3002}) and the fact that $f$ is
subadditive and satisfying (\ref{ineq}).
\end{proof}

\begin{rem}\label{add}
\rm Let $\overline{Q}_n$ be a conditional symmetrization of $Q_{n}$
given $M_{n}$. Then a tail inequality similar to (\ref{symm}) holds for
$\sup_{k\ge 1}|\Pi_{k-1}\overline{Q}_{k}|$ as well. However, in contrast to the
$\overline{Q}_{k}^{(2)}$, the $\overline{Q}_{k}$ may be degenerate in which
case an analog of (\ref{tailin}) does not follow. This is the reason for
considering $\sup_{k\ge 1}|\Pi_{2k-2}\overline{Q}_{k}^{(2)}|$ in the above
lemma.
\end{rem}

\begin{lemma}\label{e70} If $0<\mmp\{|M|<1\}\le\mmp\{|M|\le 1\}=1$, then
\begin{equation}\label{er5001}
\me\sigma(x)\ =\ 1+\sum_{n=1}^\infty \mmp\{|\Pi_n|>x\}\ \le\ 2J\big(|\log x|\big),
\end{equation}
for each $x\in (0,1]$, where $\sigma(x)\gl\inf\{n\ge 1:|\Pi_{n}|<x\}$.
Furthermore, for any $\eta>0$ such that
$$ \alpha\ \gl\ \mmp\Big\{\sup_{n\ge 1}|\Pi_{n-1}Q_n|\le\eta\Big\}\ >\ 0, $$
the function $V(x)\gl1+\sum_{n=1}^\infty \mmp\Big\{\displaystyle{
\max_{1\le k\le n}}
|\Pi_{k-1}Q_k|\le\eta, |\Pi_{n}|>x\Big\}$ satisfies
\begin{equation}\label{er5000}
V(x)\ \ge\ \alpha J\big(|\log x|\big)
\end{equation}
for each $x\in (0,1]$.
\end{lemma}
\begin{proof} Inequality (\ref{er5001}) was proved in \cite{Er73}. Below we use
the idea of an alternative proof of this result given on p.\ 153-154 in
\cite{Chow}.

Given our condition on $M$, the sequence $S_{n}=-\log|\Pi_{n}|$, $n=0,1,...$,
forms a random walk with nondegenerate increment distribution
$\mmp\{\xi\in\cdot\}$, $\xi\gl -\log |M|$. For $x>0$, put further $S_{0}^{(x)}\gl 0$
and $S_{n}^{(x)}\gl\sum_{k=1}^{n}(\xi_{k}\wedge x)$ for $n=1,2,...$, where the
$\xi_{k}$ are independent copies of $\xi$. Let
$$ T_{x}\ \gl\ \inf\Big\{n\geq 1: S_n\geq x \ \text{or} \max_{1\le k\le n}
|\Pi_{k-1}Q_k|>\eta\Big\}. $$
Then
$$\me T_{x}\ =\ \sum_{n\ge 1}\mmp\{T_x\geq n\}\ =\ V(e^{-x}) $$
and Wald's identity provide us with
\begin{equation}\label{er3}
\me S_{T_x}^{(x)}\ =\ \me (\xi\wedge x)\,\me T_x\ =\ A(x)V(e^{-x}).
\end{equation}
Putting $B\gl\{\sup_{k\ge 1}|\Pi_{k-1}Q_{k}|\le\eta\}$, we also have
$$ x\,\1_{B}\ \le\ (S_{T_x}\wedge x)\,\1_{B}\ \le\ S_{T_x}\wedge x
\ \le\ S_{T_x}^{(x)}. $$
Consequently,
$$ \me S_{T_x}^{(x)}\ \ge\ \alpha x, $$
which in combination with (\ref{er3}) implies (\ref{er5000}).
\end{proof}

\begin{lemma}\label{bas} Suppose $M,Q\ge 0$ a.s.\ and $0<\mmp\{M<1\}\le
\mmp\{M\le 1\}=1$. Let $f$ be the function defined in Lemma \ref{function}. Then
$$\me f\Big(\sup_{n\ge 1}\Pi_{n-1}Q_n\Big)<\infty\quad\Rightarrow\quad
\me f(Q)J(\log^{+}Q)<\infty. $$
\end{lemma}
\begin{proof} We first note that the moment assumption and
$\lim_{x\to\infty}f(x)=\infty$ together ensure
$\sup_{n\ge 1}\Pi_{n-1}Q_n<\infty$ a.s. Therefore, there exists an $\eta>1$ such
that $\alpha=\mmp\{\sup_{n\ge 1}\Pi_{n-1}Q_n \le\eta\}>0$.
We further point out that the monotonicity of $f$ and (\ref{ineq}) imply
$f(Q^{1/2})\ge C f(Q/2)$ for some $C\in (0,1)$.

Now fix any $\gamma>\eta$ and infer for $x\ge\eta$ (with $V$ as in the previous
lemma)
\begin{eqnarray*}
&&\mmp\Big\{\sup_{n\ge 1}\Pi_{n-1}Q_n >x\Big\}\\
&=&\mmp\{Q_1>x\}+\sum_{n\ge 1}
\mmp\Big\{\max_{1\le k\le n}\Pi_{k-1}Q_k\le x,\,\Pi_{n}Q_{n+1}>x\Big\}\\
&\ge&\mmp\{Q_1>\gamma x\}+\sum_{n\ge 1}\mmp\Big\{\max_{1\le k\le n}
\Pi_{k-1}Q_k \le\eta,\,\Pi_n Q_{n+1}>x,\,Q_{n+1}>\gamma x\Big\}\\
&\ge&\int_{\gamma x}^\infty\left(1+\sum_{n\ge 1}
\mmp\Big\{\max_{1\le k\le n}\Pi_{k-1}Q_k \le\eta,\,
\Pi_n>x/y\Big\}\right)\ \mmp\{Q\in dy\}\\
&=&\me V(x/Q)\1_{\{Q>\gamma x\}}\\
&\ge&\alpha\,\me J\big(|\log(x/Q)|\big)\1_{\{Q>\gamma x\}},
\end{eqnarray*}
the last inequality following by Lemma \ref{e70}. With this at hand, we further
obtain
\begin{eqnarray*}
\infty &>&\me f\Big(\sup_{n\ge 1}\Pi_{n-1}Q_n\Big)\\
&\ge&\int_{\eta}^{\infty}f'(x)\mmp\Big\{\sup_{n\ge 1}\Pi_{n-1}Q_n >x\Big\}\ dx\\
&\ge&\alpha\int_{\eta}^{\infty}f'(x)\,\me J\big(|\log(x/Q)|\big)
\1_{\{Q>\gamma x\}}\ dx\\
&=&\alpha\,\me\left(\int_{\eta}^{Q/\gamma}f'(x)J\big(|\log(x/Q)|\big)\ dx
\right)\\
&\ge&\alpha\,\me\left(\1_{\{Q>\gamma^{2}\}}\int_{\eta}^{Q^{1/2}}
f'(x)J\big(|\log(x/Q)|\big)\ dx\right)\\ 
&\ge&\alpha\,\me\left(\1_{\{Q>\gamma^{2}\}}f(Q^{1/2})J\bigg(\frac{\log Q}{2}
\bigg)\right)\\
&\ge&\alpha C\,\me\left(\1_{\{Q>\gamma^{2}\}}f(Q/2)J\bigg(\frac{\log Q}{2}
\bigg)\right)
\end{eqnarray*}
and this proves the assertion because $f(x)J(\log x)$ is slowly varying at
infinity by Lemma \ref{slowvar}.
\end{proof}

\begin{lemma}\label{collection} Suppose $\lim_{n\to\infty}\Pi_n=0$ a.s. Let $f$ be
the function defined in (\ref{function}), $\sigma$ the ladder epoch defined in
(\ref{ladder}) and $\sigma^{*}\gl\inf\{n\ge 1:|\Pi_{n}|>1\}$ its dual. Then
the following assertions are equivalent:
\begin{eqnarray}
&&\me f\big(|M|\big)J\big(\log^{+}|M|\big)<\infty.\label{cond1}\\
&&\me f(|\Pi_{\sigma^{*}}|)1_{\{\sigma^{*}<\infty\}}<\infty,\label{cond2}\\
&&\me f\Big(\suin|\Pi_n|\Big)<\infty,\label{cond3}\\
&&\me f\Big(\sup_{0\le n<\sigma}|\Pi_{n}|\Big)
J\Big(\sup_{0\le n<\sigma}\log^{+}|\Pi_{n}|\Big)<\infty,\label{cond4}
\end{eqnarray}
\end{lemma}
\begin{rem} \rm Rewriting Lemma \ref{collection} in terms of
$S_{n}=-\log|\Pi_{n}|$, $n=0,1,...$ and the function $b$ (recalling that
$b(\log x)\sim f(x)$), the result appears to be known under additional
restrictions on
$\{S_{n}:n=0,1,...\}$ and/or $b$, see Theorem 1 of \cite{Janson} for the
case $\me S_1\in (-\infty, 0)$ and $b$ an\nobreak (increa\-sing) power function,
Theorem 3 of
\cite{Als} for the case $\me S_1\in (-\infty, 0)$ and regularly varying $b$,
and Proposition 4.1 of \cite{KestMal} for the case $S_{n}\to -\infty$
a.s. and $b$ again a power function. In view of these results, our main
contribution is the proof of "(\ref{cond3})$\Rightarrow$(\ref{cond4})"
with the help of Lemma \ref{bas}.
\end{rem}
\begin{proof}
The equivalence "(\ref{cond1}) $\Leftrightarrow$ (\ref{cond2}) $\Leftrightarrow$
(\ref{cond3})", rewritten in terms of $\{S_{n}:n=0,1,...\}$
and $b$, takes the form
\begin{eqnarray*}
\me b\Big(\sup_{0\le n<\sigma}S_{n}\Big)J\Big(\sup_{0\le n<\sigma}S_{n}\Big)
<\infty\ &\Leftrightarrow&\
\me b(S_{\sigma^{*}})\1_{\{\sigma^{*}<\infty\}}<\infty\\
&\Leftrightarrow&\ \me b\Big(\sup_{n\ge 0}S_{n}\Big)<\infty,
\end{eqnarray*}
where $b$ is regularly varying with index $\alpha>0$. A proof for the special case
$b(x)=x^{\alpha}$ can be found in \cite{KestMal}, as mentioned above.
But the arguments given there are easily seen to hold for regularly varying $b$
as well whence further details are omitted here.

"(\ref{cond3})$\Rightarrow$(\ref{cond4})". Define the sequence $(\sigma_{n})_
{n\ge 0}$ of ladder epochs associated with $\sigma$, given by $\sigma_{0}\gl 0$,
$\sigma_{1}\gl\sigma$ and (recalling $\Pi_{k:l}=M_{k}\cdot...\cdot M_{l}$)
$$ \sigma_{n}\ \gl\ \inf\{k>\sigma_{n-1}:|\Pi_{\sigma_{n-1}:k}|\le 1\} $$
for $n\ge 2$. Put further
\begin{eqnarray*}
\widehat{\Pi}_{n}^{*}&\gl&\sup\{|\Pi_{\sigma_{n-1}}|,|\Pi_{\sigma_{n-1}+1}|,...,
|\Pi_{\sigma_{n}-1}|\},\\
\widehat{M}_{n}&\gl&\prod_{j=\sigma_{n-1}+1}^{\sigma_{n}}|M_{j}|,\\
\widehat{\Pi}_{n}&\gl&\prod_{j=1}^{n}\widehat{M}_{j}\ =\ \Pi_{\sigma_{n}}\\
\widetilde{Q}_{n}&\gl&1\vee\sup\big\{|\Pi_{\sigma_{k-1}+1:\sigma_{k-1}+k}|:1\le k
\le\sigma_{n}-\sigma_{n-1}\big\}.
\end{eqnarray*}
for $n=0,1,...$
The random vectors $(\widehat{M}_n,\widetilde{Q}_n),\ n=1,2,...$ are independent
copies of $(\widehat{M},\widetilde{Q})\gl(|\Pi_{\sigma}|,\sup_{0\le k<\sigma}
|\Pi_{k}|)$. Moreover, $\widehat{\Pi}_{n}^{*}=|\Pi_{\sigma_{n-1}}|
\widetilde{Q}_{n}=\widehat{\Pi}_{n-1}\widetilde{Q}_{n}$ and
$$ \sup_{n\ge 0}|\Pi_n|\ =\ \sup_{n\ge 1}|\widehat{\Pi}_{n}^{*}|
\ =\ \sup_{n\ge 1}\widehat{\Pi}_{n-1}\widetilde{Q}_{n}. $$
As, by construction, $\mmp\{\widehat{M}\le 1\}=1$ and
$\mmp\{\widehat{M}=1\}=0$, Lemma \ref{bas} enables us to conclude that
$\me f(\sup_{n\ge 0}|\Pi_n|)=\me f(\sup_{n\ge 1}
\widehat{\Pi}_{n-1}\widetilde{Q}_n)<\infty$ implies
$\me f\big(\widetilde{Q}\big)J\big(\log^{+}\widetilde{Q}\big)<\infty$ which is the
desired result.

Finally, "(\ref{cond4})$\Rightarrow$(\ref{cond1})" follows from the obvious
inequality $\sup_{0\le n<\sigma}|\Pi_{n}|$ $\ge |M_{1}|\vee 1$ and the fact that
$f(x)J(\log x)$ is nondecreasing.
\end{proof}

\section{Proof of Theorem \ref{perp}.}\label{pro1}

\emph{Sufficiency}. As condition (\ref{Q}) clearly implies $\me
J\big(\log^{+}|Q|\big)<\infty$ we infer
$\zi^{*}<\infty$ a.s.\ from Proposition \ref{exper}. Notice that our given
assumption $\lim_{n\to\infty}\Pi_n=0$ a.s.\ is valid if, and only if, one of the
following cases holds true:
\begin{itemize}
\item[(C1)] $\mmp\{|M|\le 1\}=1$ and $\mmp\{|M|<1\}>0$.
\item[(C2)] $\mmp\{|M|>1\}>0$ and $\lim_{n\to\infty}\Pi_n=0$ a.s.
\end{itemize}
We will consider these cases separately, in fact Case (C2) will be handled by
reducing it to the first case via an appropriate stopping argument.
\bigskip

\textbf{Case (C1)}: We will prove (\ref{Z*}) or, equivalently, $\me
f(\zi^{*})<\infty$. According to Lemma \ref{function}, (\ref{Q}) is
equivalent to
\begin{equation}\label{Q100}
\me f(|Q|)J\big(\log^{+}|Q|\big)<\infty
\end{equation}
which in view of (P2) particularly ensures $\me f(|Q|)<\infty$.

Using the properties of $f$ stated in Lemma \ref{function} (which particularly
ensure subadditivity) and $\sup_{n\ge 0}|\Pi_{n}|=|\Pi_{0}|=1$, we
obtain for fixed $a\in (0,1)$
\begin{eqnarray*}
\me f(\zi^{*})&=&\lim_{n\to\infty}\me f\left(
\sum_{k=1}^{n}|\Pi_{k-1}Q_{k}|\right)\\
&\le&\lim_{n\to\infty}\sum_{k=1}^{n}\me f(|\Pi_{k-1}Q_{k}|)\\
&\le&\int_{0}^{\infty}f'(x)\sum_{k\ge 1}
\mmp\{|\Pi_{k-1}Q_k|>x\}\ dx\\
&=&\int_{0}^{\infty}f'(x)\sum_{k\ge 1}
\mmp\{|\Pi_{k-1}Q_k|>x,|Q_k|>x/a\}\ dx\\
&+&\int_{0}^{\infty}f'(x)\sum_{k\ge 1}
\mmp\{|\Pi_{k-1}Q_k|>x,x<|Q_k|\le x/a\}\ dx\\
&=&I_{1}+I_{2}
\end{eqnarray*}
The second integral is easily estimated with the help of (\ref{er5001}) as
\begin{eqnarray*}
I_{2}&\le&\left(\sum_{k\ge 1}\mmp\{|\Pi_{k-1}|>a\}\right)\int_{0}^{\infty}
f'(x)\,\mmp\{|Q|>x\}\ dx\\
&\le&2J(|\log a|)\,\me f(|Q|)\ <\ \infty,
\end{eqnarray*}
so that we are left with an estimation of $I_{1}$.

The concavity of $f$ in
combination with $f(0)=0$ and $f'(0+)<\infty$ (see Lemma \ref{function}) gives
$f(x)\le f'(0+)x$ for all $x>0$. As in Lemma \ref{e70}, let $\sigma(t)=\inf\{n\ge
1:|\Pi_{n}|<t\}$ for $t>0$ and recall from there that $\me\sigma(t)\le 2J(|\log t|)$ for
$t\le 1$. For $t>1$, we trivially have $\sigma(t)\equiv 1$. Finally, put
$\rho\gl\me|M|$, so that $\rho\in (0,1)$ and furthermore $\sum_{k\ge 1}\me|\Pi_{k}|
=(1-\rho)^{-1}$. Hence
$$ \sum_{k\ge 1}\me f(|\Pi_{k}|)\ \le\ \Lambda\ \gl\ {f'(0+)\over 1-\rho}
\ <\ \infty. $$
By combining these facts, we infer
\begin{eqnarray*}
I_{1}&=&\int_{0}^{\infty}f'(x)\int_{(x/a,\infty)}\sum_{k\ge 1}
\mmp\{|\Pi_{k-1}|>x/y\}\ \mmp\{|Q|\in dy\}\ dx\\
&=&\int_{(0,\infty)}\int_{0}^{a}yf'(xy)\sum_{k\ge 0}
\mmp\{|\Pi_{k}|>x\}\ dx\ \mmp\{|Q|\in dy\}\\
&\le&\int_{(0,\infty)}\sum_{k\ge 0}\me f\big(y(|\Pi_{k}|\wedge a)\big)
\ \mmp\{|Q|\in dy\}\\
&\le&\int_{(1,\infty)}\sum_{k\ge 0}\me f\big(y(|\Pi_{k}|)\big)
\ \mmp\{|Q|\in dy\}\ +\ \sum_{k\ge 0}\me f(|\Pi_{k}|)\\ \noalign{\break}
&\le&\int_{(1,\infty)}\left[f(y)\,\me\sigma(1/y)+\me\Bigg(\sum_{k\ge\sigma(1/y)}
f(y|\Pi_{k}|)\Bigg)\right]\ \mmp\{|Q|\in dy\}\ +\ \Lambda\\
&\le&\int_{(1,\infty)}\left[f(y)\,\me\sigma(1/y)+\me\Bigg(\sum_{k\ge\sigma(1/y)}
f(|\Pi_{\sigma(1/y)+1:k}|)\Bigg)\right]\ \mmp\{|Q|\in dy\}\ +\ \Lambda\\
&=&\int_{(1,\infty)}\left[f(y)\,\me\sigma(1/y)+\me\left(\sum_{k\ge 0}
f(|\Pi_{k}|)\right)\right]\ \mmp\{|Q|\in dy\}\ +\ \Lambda\\
&\le&\int_{(1,\infty)}2f(y)J(|\log y|)\ \mmp\{|Q|\in dy\}\ +\ 2\Lambda\\
&\le&2\,\me f(|Q|)J(\log^{+}|Q|)\ +\ 2\Lambda.
\end{eqnarray*}
But the final line is clearly finite by our given moment assumptions
which completes the proof for Case (C1).
\medskip

\textbf{Case (C2)}: As already announced, we will handle this case by using
a stopping argument based on the ladder epoch $\sigma$ given in (\ref{ladder}).
We adopt the notation of the proof of Lemma \ref{collection},
in particular $(\sigma_{n})_{n\ge 0}$ denotes the sequence of successive ladder
epochs associated with $\sigma$. Put further
$$ \widehat{Q}_{n}\ \gl\ \sum_{k=\sigma_{n-1}+1}^{\sigma_{n}}|\Pi_{\sigma_{n-1}
+1:k-1}Q_{k}| $$
for $n\ge 1$ which are independent copies of $\widehat{Q} \gl\widehat{Q}_{1}=
Q^{(\sigma)}$. Notice that
\begin{equation}\label{ladperp}
\zi^{*}\ =\ \sum_{k\ge 1}\widehat{\Pi}_{k-1}\widehat{Q}_{k}.
\end{equation}
It will be shown now that condition (\ref{Q100}) holds true with $\widehat{Q}$
instead of $Q$. Since
$\widehat{M}=|\Pi_{\sigma}|\in (0,1)$ a.s.\ and thus satisfies the condition of Case
(C1), we then arrive at the desired conclusion $\me f(\zi^{*})<\infty$.

By Lemma \ref{subadd}, there is a subadditive $\phi(x)$ of the same asymptotic
behavior as $f(x)J(\log x)$, as $x\to\infty$. Hence it suffices to verify
$\me\phi(\widehat{Q})<\infty$. Use the obvious inequality
$$ \widehat{Q}\ \le\ \sup_{1\le k\le\sigma}|\Pi_{k-1}|\sum_{k=1}^{\sigma}
|Q_{k}|\ =\ \widetilde{Q}\sum_{k=1}^{\sigma}|Q_{k}|. $$
in combination with property (\ref{submulti}) and the subadditivity of $\phi$ to
infer
$$ \me\phi(\widehat{Q})\ \le\ C\left(\me\phi(\widetilde{Q})+\me\left(
\sum_{k=1}^{\sigma}\phi(|Q_{k}|)\right)\right). $$
But the right hand expression is finite because $\me\phi(\widetilde{Q})<\infty$
is ensured by (\ref{M}) and Lemma \ref{collection} and because
$$ \me\left(\sum_{k=1}^{\sigma}\phi(|Q_{k}|)\right)\ =\ \me\phi(|Q|)\,
\me\sigma\ <\ \infty $$
follows from Wald's identity, condition (\ref{Q}) and $\me\sigma<\infty$ which in
turn is a consequence of our assumption $\lim_{n\to\infty}\Pi_{n}=0$ a.s.
\medskip

\noindent
\emph{Necessity}. This is easier. Assuming (\ref{Z}) or, equivalently,
$\me f(|\zi|)<\infty$, we infer from Lemma \ref{perpsup}
$$ \me f\Big(\sup_{n\ge 1}|\widetilde{\Pi}_{n-1}Q_{n}|\Big)\ \le
\ \me f\Big(\sup_{n\ge 1}|\Pi_{n-1}Q_{n}|\Big)\ <\ \infty, $$
where $\widetilde{\Pi}_{n}\gl\prod_{k=1}^{n}(M_{k}\wedge 1)$, and thereupon
$\me f(|Q|)J(\log^{+}|Q|)<\infty$ by Lemma \ref{bas} (as
$\mmp\{|M\wedge 1|<1\}=\mmp\{|M|<1\}>0$).

Left with the proof of (\ref{M}), we get $\me f(\sup_{n\ge 0}|\Pi_{n}|)<\infty$
by another appeal to Lemma \ref{perpsup} and then the assertion by invoking
Lemma \ref{collection}. This completes the proof of Theorem \ref{perp}.
\hfill$\square$

\section{Size-biasing and the results for $W_n$}\label{ui}

\subsection{Modified branching random walk}\label{spine1}

We adopt the situation described in Subsection \ref{brw}.
Recall that $\mm$ denotes a generic copy of the point process describing
the displacements of children relative to its mother in the considered population.
The following construction of the associated {\it modified BRW} with a
distinguished ray $(v_{n})_{n\ge 0}$, called {\it spine}, is based on
\cite{BiggKypr04} and \cite{Lyons}.

Let $\mm^\ast$ be
a point process whose law has Radon-Nikodym derivative $m(\gamma)^{-1}$
$\sum_{i=1}e^{\gamma X_{i}}$ with
respect to the law of $\mm$. The individual $v_0=\varnothing$ residing at the
origin of the real line has children, the displacements of which relative to
$v_0$ are given by a copy $\mm_{0}^{*}$ of $\mm^\ast$. All the children of $v_0$
form the first generation of the population, and among these the spinal successor
$v_{1}$ is picked with a probability proportional to $e^{\gamma s}$ if $s$ is the
position of $v_{1}$ relative to $v_{0}$ (size-biased selection). Now, while $v_1$
has children the displacements of which relative to $v_1$ are given by another
independent copy $\mm_{1}^{*}$ of $\mm^\ast$, all other individuals of the first
generation produce and spread offspring according to independent copies of $\mm$
(i.e., in the same way as in the given BRW). All children of the individuals of the
first generation form the second generation of the population, and among the
children of $v_{1}$ the next spinal individual $v_{2}$ is picked with probability
$e^{\gamma s}$ if $s$ is the position of $v_{2}$ relative to $v_{1}$. It produces
and spreads offspring according to an independent copy $\mm_{2}^{*}$ of $\mm^{*}$
whereas all siblings of $v_{2}$ do so according to independent copies of $\mm$,
and so on. Let
$\widehat{\mm}_{n}$ denote the point process describing the positions of all
members of the $n$-th generation. We call
$\{\widehat{\mm}_{n}:n=0,1,...\}$ a {\it modified BRW} associated with
the ordinary BRW $\{\mm_{n}:n=0,1,...\}$.

Recall that $\mathbf{T}$ denotes the Galton-Watson tree associated with
$\{\mm_{n}:n=0,1,...\}$, and denote by $\whbfT$ the corresponding
size-biased tree associated with $\{\widehat{\mm}_{n}:n=0,1,...\}$. Let $\bfP$ be
the distribution of the random weighted tree $(\mathbf{T},\mathbf{L})$, where
$\mathbf{L}
\gl (L(v))_{v\in\mathbf{T}}$ with $L(v)=e^{\gamma S(v)}/m(\gamma)^{|v|}$ denoting
the weight (as defined in Subsection \ref{brw}) attached to the node $v$ residing
at $S(v)$. Similarly, let $\whL(v)
\gl e^{\gamma \whS(v)}/m(\gamma)^{|v|}$ be the weight of any $v\in
\whbfT$
if $\whS(v)$ denotes its position, i.e., $\widehat{\mm}_{n}=\sum_{v\in\mathbf{T}
^{*}:|v|=n}\delta_{\whS(v)}$ for each $n=0,1,...$ The distribution of the thus
obtained random weighted tree $(\whbfT,\whbfL)$,
$\whbfL\gl (\whL(v))_{v\in\whbfT}$,
is denoted as $\bfQ$. Both, $\bfP$ and $\bfQ$, are probability measures on
the space
$$ {\Bbb W}\ \gl\ \{(t,l):t\subset\mathbf{V}\} $$
of weighted subtrees of $\mathbf{V}$ with the same root,
where $l:t\to\mr$ is the weight function putting weight $l(v)$ to each $v\in t$.
Endow this space with the filtration $\{\mathcal{G}_{n}:n=0,1,...\}$, where
$\mathcal{G} _{n}$ is generated by the sets
$$ [t,l]_{n}\ \gl\ \{(t',l')\in {\Bbb W}:t_{n}=t_{n}'\hbox{ and }l_{|t_{n}}
=l'_{|t_{n}}\},\quad (t,l)\in {\Bbb W}. $$
Here $t_{n}\gl\{v\in t:|v|\le n\}$. Put further $\mathcal{G}\gl\sigma\{
\mathcal{G}_{n}:n=0,1,...\}$. Then the mappings $z_{n},w_{n}:{\Bbb W}\to
[0,\infty)$, defined as
$$ z_{n}(t,l)\ \gl\ \sum_{v\in t_{n}}l(v)\quad\hbox{and}\quad
w_{n}(t,l)\ \gl\ m(\gamma)^{-n}z_{n}(t,l), $$
are $\mathcal{G}_{n}$-measurable for each $n\ge 0$, and we have
$$ W_{n}\ =\ w_{n}\circ (\mathbf{T},{\bf L}),\quad n=0,1,... $$
Put also $\whW_{n}\gl w_{n}\circ (\whbfT,\whbfL)$ and
$\whW\gl\limsup_{n\to\infty}\whW_{n}$. Then
\begin{eqnarray}
&&\mmp((W_{n})_{n\ge 0}\in\cdot)=\bfP((w_{n})_{n\ge 0}\in\cdot)\label{d1}\\
&\hbox{and}&\mmp((\whW_{n})_{n\ge 0}\in\cdot)=\bfQ((w_{n})_{n\ge 0}
\in\cdot).\label{d2}
\end{eqnarray}

The relevance of these definitions with view to the
martingale $\{W_{n}:n=0,1,...\}$ to be studied hereafter is provided by the
following lemma (see Prop.\ 12.1 and Thm.\ 12.1 in \cite{BiggKypr04}).

\begin{lemma}\label{sizeb} For each $n\ge 0$, $w_{n}$ is the Radon-Nikodym
derivative of $\bfQ$ with respect to $\bfP$ on $\mathcal{G}_n$. Moreover,
if $w\gl\limsup_{n\to\infty}w_{n}$, then
\item{(1)\quad} $w_{n}$ is a $\bfP$-martingale and
$1/w_{n}$ is a $\bfQ$-martingale.
\item{(2)\quad} $\me W=\me_{\bfP}w=1$ if and only if $\mmp\{\whW<\infty\}=
\bfQ\{w<\infty\}=1$.
\item{(3)\quad} $\me W=\me_{\bfP}w=0$ if and only if $\mmp\{\whW=\infty\}=
\bfQ\{w=\infty\}=1$.
\end{lemma}

Let us point out that, in view of (\ref{d1}) and (\ref{d2}), the first assertion
of Lemma \ref{sizeb}(1) just restates the martingale property of $W_{n}$, while
the second one says that the same holds true for $1/\whW_{n}$.
The link between $W_{n}$ and $\whW_{n}$ is provided by

\begin{lemma}\label{link1}
For each $n=0,1,...$, $\whW_{n}$ is a size-biasing of $W_{n}$, that is
\begin{equation}\label{link1.1}
\me W_{n}f(W_{n})\ =\ \me f(\whW_{n}).
\end{equation}
for each function $f:\mr^{+}\to\mr^{+}$. More generally,
\begin{equation}\label{link1.2}
\me W_{n}h(W_{0},...,W_{n})\ =\ \me h(\whW_{0},...,\whW_{n}).
\end{equation}
for each Borel function $h:(\mr^{+})^{n+1}\to\mr^{+}$.
\end{lemma}

\begin{proof}
It suffices to note that, by Lemma \ref{sizeb},
\begin{eqnarray*}
\me W_{n}h(W_{0},...,W_{n})&=&\me_{\bfP}w_{n}h(w_{0},...,w_{n})\\
&=&\me_{\bfQ}h(w_{0},...,w_{n})
\ =\ \me h(\whW_{0},...,\whW_{n})
\end{eqnarray*}
for each $n=0,1,...$ and $h$ as stated, where the $\mathcal{G}_{n}$-measurability
of $(w_{0},...,w_{n})$ should be observed for the second equality.
\end{proof}

\subsection{Connection with perpetuities}\label{connex}

Next we have to make the connection with perpetuities. For $u\in\whbfT$,
let $\whcalN(u)$ denote the set of children of $u$ and, if $|u|=k$,
$$ \whW_{n}(u)\ =\ \sum_{v:uv\in\whbfT_{k+n}}{\whL(uv)\over
\whL(u)},\quad n=0,1,... $$
Since all individuals off the spine reproduce and spread as in the unmodified
BRW, we have that the $\{\whW_{n}(u):n=0,1,...\}$ for $u\in\bigcup_{n\ge 0}
\whcalN(v_{n})\backslash\{v_{n+1}\}$ are independent copies of
$\{W_{n}:n=0,1,...\}$. For $n\in\mn$, define further
\begin{equation}\label{eqM}
M_{n}\ \gl\ {\whL(v_{n})\over \whL(v_{n-1})}\ =\ {e^{\gamma(\whS(v_{n})
-\whS(v_{n-1}))}\over m(\gamma)}
\end{equation}
and
\begin{equation}\label{eqQ}
Q_{n}\ \gl\ \sum_{u\in\whcalN(v_{n-1})}
{\whL(u)\over \whL(v_{n-1})}\ =\ \sum_{u\in\whcalN(v_{n-1})}
{e^{\gamma(\whS(u)-\whS(v_{n-1}))}\over m(\gamma)}.
\end{equation}
Then it is easily checked that the $\{(M_{n},Q_{n}):n=1,2,...\}$ are i.i.d.\
with distribution given by
\begin{eqnarray*}
\mmp\{(M,Q)\in A\}&=&\me\Bigg(\sum_{i=1}^{N}{e^{\gamma X_{i}}\over m(\gamma)}
\1_{A}\Bigg({e^{\gamma X_{i}}\over m(\gamma)},\sum_{j=1}^{N}{e^{\gamma X_{j}}
\over m(\gamma)}\Bigg)\Bigg)\\
&=&\me\Bigg(\sum_{|u|=1}L(u)\1_{A}\Bigg(L(u),\sum_{|v|=1}L(v)\Bigg)\Bigg)
\end{eqnarray*}
for any Borel set $A$, where $(M,Q)$ denotes a generic copy of $(M_{n},
Q_{n})$ and our convention $\sum_{|u|=n}=\sum_{u\in\mathbf{T}_{n}}$ should be
recalled from Section 1. In particular,
\begin{equation*}
\mmp\{Q\in B\}\ =\ \me\Bigg(\sum_{|u|=1}L(u)\1_{B}\Bigg(
\sum_{|u|=1}L(u)\Bigg)\Bigg)\ =\ \me W_{1}\1_{B}(W_{1})
\end{equation*}
for any measurable $B$, that is
\begin{equation}\label{qw1}
\mmp\{Q\in dx\}\ =\ x\,\mmp\{W_{1}\in dx\}.
\end{equation}
Notice that this implies
\begin{equation}\label{positive}
\mmp\{Q=0\}=0.
\end{equation}
As for the distribution of $M$, we have
\begin{equation*}
\mmp\{M\in B\}\ =\ \me\Bigg(\sum_{|u|=1}L(u)\1_{B}(L(u))\Bigg)
\end{equation*}
which is in accordance with the definition given in (\ref{Zm}). As we see
from (\ref{eqM}),
\begin{equation}\label{eqPi}
\Pi_{n}\ =\ M_{1}\cdot...\cdot M_{n}\ =\ \whL(v_{n}),\quad n=0,1,...
\end{equation}

Here is the lemma that provides the connection between the sequence
$\{\whW_{n}:n=0,1,...\}$ and the perpetuity generated by $\{(M_{n},Q_{n}):
n=0,1,...\}$. Let $\mathcal{A}$ be the $\sigma$-field generated by
$\{(M_{n},Q_{n}):n=0,1,...\}$ and the family $\{\mm_{n}^{*}:
n=0,1,...\}$, where $\mm_{n}^{*}$ is the copy of $\mm^{*}$ describing the
displacement of the children of $v_{n}$ relative to its mother. For $n\ge 1$ and
$k=1,...,n$, put also
$$ R_{n,k}\ \gl\ \sum_{u\in\whcalN(v_{k-1})\backslash\{v_{k}\}}
{\whL(u)\over\whL(v_{k-1})}\Big(\whW_{n-k}(u)-1\Big) $$
and notice that $\me\big(R_{n,k}|\mathcal{A}\big)=0$ because each $\whW_{n-k}(u)$
is independent of $\mathcal{A}$ with mean one.

\begin{lemma}\label{link2}
With the previous notation the following identities hold true for each $n\ge 0$:
\begin{equation}\label{link2.1}
\whW_{n}\ =\ 1+\sum_{k=1}^{n}\Pi_{k-1}\big(Q_{k}+R_{n,k}\big)
\end{equation}
and
\begin{equation}\label{link2.2}
\me\big(\whW_{n}|\mathcal{A}\big)\ =\ 1+\sum_{k=1}^{n}\Pi_{k-1}Q_{k}\quad
\mmp\hbox{-a.s.}
\end{equation}
\end{lemma}

\begin{proof}
Each $v\in\whbfT_{n}$ has a most recent ancestor in $\{v_{k}:k=0,1,...\}$. By
using this and recalling (\ref{eqQ}) and (\ref{eqPi}), one can easily see that
\begin{eqnarray*}
\whW_{n}&=&\whL(v_{n})+\sum_{k=1}^{n}\sum_{u\in\whcalN(v_{k-1})
\backslash\{v_{k}\}}\whL(u)\whW_{n-k}(u)\\
&=&\Pi_{n}+\sum_{k=1}^{n}\Pi_{k-1}\bigg(Q_{k}-{\whL(v_{k})\over
\whL(v_{k-1})}+1+R_{n,k}\bigg)\\
&=&\Pi_{n}+\sum_{k=1}^{n}\big(\Pi_{k-1}-\Pi_{k}\big)+\sum_{k=1}^{n}\Pi_{k-1}
\big(Q_{k}+R_{n,k}\big)
\end{eqnarray*}
which obviously gives (\ref{link2.1}) as $\Pi_{0}=1$.
But the second assertion is now immediate when observing that $\me(\Pi_{k-1}R_{n,k}|
\mathcal{A})=\Pi_{k-1}\me(R_{n,k}|\mathcal{A})=0$ a.s.
\end{proof}

\subsection{Two further auxiliary results}

We continue with two further auxiliary results about the martingale $W_{n}$
and its size-biasing $\whW_{n}$.

\begin{lemma}\label{count} Let $W^{*}\gl\sup_{n\ge 0}W_{n}$ and
$\whW^{*}\gl\sup_{n\ge 0}\whW_{n}$. Then, for each $a\in (0,1)$, there exists
$b=b(a)\in\mr^+$ such that for all $t>1$
\begin{equation}\label{tailsup}
\mmp\{W>t\}\ \le\ \mmp\{W^{*}>t\}\ \le\ b\,\mmp\{W>at\}.
\end{equation}
As a consequence
\begin{equation*}
\me f(W)<\infty\quad\Leftrightarrow\quad\me f(W^{*})<\infty
\end{equation*}
for any non-negative nondecreasing concave function $f$. Replacing
$(W,W^{*})$\nobreak\ with $(\whW,\whW^{*})$, the same conslusions hold true (with $b/a$
instead of
$b$).
\end{lemma}

\begin{proof} (\ref{tailsup}) is due to Biggins \cite{Bigg} for the case of a.s.\ finite
branching (see Lemma 2 there) and has been obtained without this restriction as Lemma 1
in \cite{IksNegad} by a different argument. Its counterpart for $(\whW,\whW^{*})$ can
be found as Lemma 3 in \cite{IksRos}, but the following argument (for the nontrivial
part) using (\ref{tailsup}) and Lemma \ref{link1} is more natural and much shorter:
\begin{eqnarray*}
\mmp\{\whW^{*}>t\}&=&\sum_{n\ge 1}\mmp\Big\{\whW_{n}=\max_{0\le k\le n}
\whW_{k},\whW_{n}>t\Big\}\\
&=&\sum_{n\ge 1}\int_{\{W_{n}=\max_{0\le k\le n}W_{k},W_{n}>t\}}
W_{n}\ d\mmp\\
&=&\me W^{*}\1_{\{W^{*}>t\}}\hskip 4cm[\hbox{by (\ref{link1.2})}]\\
&\le&\int_{0}^{\infty}\mmp\{W^{*}>x\vee t\}\ dx\\
&\le&\int_{0}^{\infty}b\,\mmp\{W>a(x\vee t)\}\ dx\\
&=&b\,\me\bigg({W\over a}\1_{\{W/a>t\}}\bigg)\\
&=&{b\over a}\,\mmp\{\whW>at\}
\end{eqnarray*}
for all $t>1$.
\end{proof}

\begin{lemma}\label{degen1} Suppose that $\{W_n:n=0,1,...\}$ is uniformly
integrable. Then the following assertions hold true:
\item{(1)} If $W_{1}=1$ a.s., then $W=\whW=1$ a.s.
\item{(2)} If $\mmp\{W_{1}=1\}<1$, then $W,\whW$ are both nondegenerate.
\end{lemma}

\begin{proof} The first statement follows, as $W_{1}=1$ a.s.\ implies
the same for each $W_{n}$, $n\ge 2$ (use $W_{n}=\sum_{|v|=n-1}L(v)W_{1}(v)$
with independent $W_{1}(v)$ which are copies of $W_{1}$ and independent of
the $L(u)$, $|u|=n-1$). Conversely, if $W$ (and thus also
$\whW$ as its size-biasing) is degenerate, then the fixed point equation
(\ref{fixp}) for
$n=1$ combined with $\me W=1$ yields
$$ 1\ =\ W\ =\ \sum_{|v|=1}L(v)W(v)\ =\ \sum_{|v|=1}L(v)\ =\ W_{1}\quad
\hbox{a.s.} $$
which completes the proof.
\end{proof}

\subsection{Proof of Theorem \ref{brwL1}}

\emph{Sufficiency}. Suppose first that (\ref{Z3}) and (\ref{mea}) hold true
which, by Proposition \ref{exper}, ensures $\sum_{k\ge 1}\Pi_{k-1}Q_{k}<
\infty$ a.s.
Since $W_{n}$ is nonnegative and a.s.\ convergent to $W$, the uniform
integrability follows if we can show $\me W=1$ or, equivalently (by Lemma
\ref{sizeb}), $\mmp\{\whW<\infty\}=\bfQ\{w<\infty\}=1$.
To this end note that,
by (\ref{link2.2}) and Fatou's lemma,
$$ \me(\liminf_{n\to\infty}\whW_{n}|\mathcal{A})\ \le\ \sum_{k\ge 1}\Pi_{k-1}
Q_{k}\ <\ \infty\quad\hbox{a.s.} $$
and thus $\liminf_{n\to\infty}\whW_{n}<\infty$ a.s. As $\{1/\whW_{n}:n=0,1,...\}$
constitutes a positive and thus a.s.\ convergent martingale (see after Lemma
\ref{sizeb}), we further infer $\whW=\liminf_{n\to\infty}\whW_{n}$ and thereupon
the desired $\mmp\{\whW<\infty\}=1$.
\medskip

\noindent
\emph{Necessity}. Assume now that $\{W_{n}:n=0,1,...\}$ is uniformly integrable,
so that $\me W=1$ and thus $\whW<\infty$ a.s.\ by Lemma \ref{sizeb}(2).
Furthermore, $\whW^{*}<\infty$ a.s.\ by Lemma \ref{count}. The inequality
\begin{equation}\label{lowerb}
\whW_{n}\ \ge\ \whL(v_{n-1})\sum_{v\in\whcalN(v_{n-1})}{\whL(v)\over\whL(v_{n-1})}
\ =\ \Pi_{n-1}Q_{n}
\end{equation}
then shows that
\begin{equation}\label{limsup}
\sup_{n\ge 1}\Pi_{n-1}Q_{n}\ \le\ \whW^{*}\ <\ \infty\quad\hbox{a.s.}
\end{equation}
which in combination with $\mmp\{M=1\}<1$ (see (\ref{Z1})) allows us to appeal
to Theorem 2.1 in \cite{GolMal} to conclude validity of (\ref{Z3}) and (\ref{mea}).
\hfill $\square$

\begin{rem}\label{limit} \rm With view to the subsequent proof of Theorem
\ref{brw_conc2} it is useful to point out that the previous proof has shown that,
if $\{W_{n}:n=0,1,...\}$ is uniformly integrable,
$\whW=\lim_{n\to\infty}\whW_{n}<\infty$ a.s.\ and
\begin{equation*}
\me(\whW|\mathcal{A})\ \le\ \zi\ \gl\ \sum_{k\ge 1}\Pi_{k-1}Q_{k}\quad\hbox{a.s.}
\end{equation*}
Consequently, if $f:\mr^{+}\to\mr^{+}$ denotes any nondecreasing and concave function,
then an application of Jensen's inequality (for conditional expectations) in
combination with (\ref{link1.1}) gives
\begin{equation}\label{link5.4}
\me Wf(W)\ =\ \me f(\whW)\ \le\ \me f(\zi).
\end{equation}
\end{rem}

\subsection{Proof of Theorem \ref{brw_conc2}}

\noindent
\emph{Sufficiency}. Let $\zi$ be defined as usual with $M_{k}$ and $Q_{k}$
as in (\ref{eqM}) and (\ref{eqQ}), respectively.
Notice that $\zi^{*}=\sum_{k\ge 1}|\Pi_{k-1}Q_{k}|=\zi$ in the present context. By
Lemma \ref{connex} and (\ref{qw1}), condition (\ref{ZZZ}) translates to
$$ \me b(\log^{+}\whW_{1})J(\log^{+}\whW_{1})\ =\ \me b(\log^{+}Q)J(\log^{+}
Q)\ <\ \infty, $$
and we may naturally replace $b(\log^{+}x)$ with the asymptotically equivalent
concave function $f$ from Lemma \ref{function}. Since
$$ M_{1}\ =\ \whL(v_{1})\ \le\ \sum_{v\in\whbfT_{1}}\whL(v)\ =\ \whW_{1}, $$
we also infer $\me f(M)J(\log^{+}M)<\infty$. Hence the desired conclusion
(\ref{ZZZZ}), equivalently $\me Wf(W)<\infty$, follows by an appeal to
Theorem \ref{perp} and (\ref{link5.4}).

\medskip\noindent
\emph{Necessity}. Suppose now uniform integrability of the $W_{n}$,
$\mmp\{W_{1}=1\}<1$ and $\me Wf(W)<\infty$ with $f$ as before. Then
$\whW<\infty$ a.s.\ and $\me f(\whW)<\infty$ by another appeal to (\ref{link1.1}).
Next, Lemma \ref{count} gives $\me f(\whW^{*})<\infty$ and  then in combination
with (\ref{lowerb})
$$ \me f\Big(\sup_{k\ge 1}\widetilde{\Pi}_{k-1}Q_{k}\Big)\ \le
\ \me f\Big(\sup_{k\ge 1}\Pi_{k-1}Q_{k}\Big)\ \le\ \me f(\whW^{*})
\ <\ \infty, $$
where $\widetilde{\Pi}_{k}=\prod_{j=1}^{k}(M_{j}\wedge 1)$ is defined as in the
proof of Theorem \ref{brwL1}, from which we further see that the uniform
integrability of the $W_{n}$ ensures $\lim_{n\to\i}\Pi_{n}=0$ a.s.\ (Theorem
\ref{brwL1}) and thus $\mmp\{0<M<1\}>0$. Consequently, we can finally invoke Lemma
\ref{bas} in combination with (\ref{qw1}) to conclude
$$ \me f(Q)J(\log^{+}Q)\ =\ \me f(\whW_{1})J(\log^{+}\whW_{1})\ =\ \me W_{1}
f(W_{1})J(\log^{+}W_{1})\ <\ \infty $$
which proves (\ref{ZZZ}).\hfill $\square$

\bigskip

\textbf{Acknowledgment.} A part of this work was done while A.
Iksanov was visiting M\"{u}nster in October/November 2006.
Grateful acknowledgment is made for financial support and
hospitality. The research of A. Iksanov was partly supported by
the DFG grant, project no.436UKR 113/93/0-1.

\end{document}